\documentstyle[12pt]{article}

\begin{document}
\author{S. V. L\"udkovsky }
\title{Algebras of non-Archimedean measures on groups}
\date{7 May 2004}
\maketitle

\par Quasi-invariant measures with values in non-Archimedean fields
on a group of diffeomorphisms were constructed for non-Archimedean
manifolds $M$ in \cite{lu1,lutmf99}.
On non-Archimedean loop groups and semigroups they were provided in 
\cite{luumn983,lubp2,luseamb2}.
A Banach space over a local field also serves as the additive group
and quasi-invariant measures on it were studied in
\cite{luumn582,lulapm}.
\par This article is devoted to the investigation of properties of
quasi-invariant measures with values in non-Archimedean fields
that are important for analysis on topological
groups and for construction of irreducible representations.
The following properties are investigated: 
\par (1) convolutions of measures and functions,
\par (2) continuity of functions of measures,
\par (3) non-associative algebras generated with the help 
of quasi-invariant measures.
The theorems given below show that many differences appear to be
between locally compact and non-locally compact groups.
The groups considered below are supposed to have structure 
of Banach manifolds over the corresponding fields. \\
\par {\bf 1. Definitions. (a).} Let $G$ be a Hausdorff separable
topological group.
A tight measure $\mu $ on $Af (G,\mu )$ with values in a non-Archimedean
field $\bf F$ is called left-quasi-invariant (or right)
relative to a dense subgroup $H$ of $G$, if $\mu _{\phi }(*)$
(or $\mu ^{\phi }(*)$) is equivalent to $\mu (*)$ 
for each $\phi \in H$, where $Bco (G)$ is the algebra of all
clopen subsets of $G$, $Af (G,\mu )$ denotes its completion by $\mu $, 
$\mu _{\phi }(A):=\mu (\phi ^{-1}A)$, $\mu ^{\phi }(A):=\mu (A\phi ^{-1})$
for each $A\in Af(G,\mu )$, $\rho _{\mu }(\phi ,g):=
\mu _{\phi }(dg)/\mu (dg) \in L(G, Af(G,\mu ), \mu , {\bf F}) $
(or ${\tilde {\rho }}_{\mu }(\phi ,g):=\mu ^{\phi }(dg)/\mu (dg)$)
denotes a left (or right) quasi-invariance factor, $\bf F$
is a non-Archimedean field complete relative to its uniformity and
such that ${\bf K}_s\subset \bf F$. We assume
that a uniformity $\tau _G$ on $G$ is such that  $\tau _G|H\subset 
\tau _H$, $(G,\tau _G)$ and $(H,\tau _H)$ are complete. We suppose
also that there exists an open base in $e\in H$ such that their closures
in $G$ are compact (such pairs exist for loop groups and
groups of diffeomorphisms and Banach-Lie groups). 
We denote by $M_l(G,H)$ (or $M_r(G,H)$) a set of left-( or right)
quasi-invariant tight measures on $G$ relative to $H$ with a finite norm
$\| \mu \| <\infty $.
\par {\bf (b).} Let $L_H(G,\mu ,{\bf F})$ denotes the Banach space
of functions $f: G\to \bf F$ such that $f_h(g)\in L(G,\mu ,{\bf F})$
for each $h\in H$ and 
$$\| f\|_{L_H(G,\mu ,{\bf F})}:=\sup_{h\in H}
\| f_h\|_{L(G,\mu ,{\bf F})}<\infty ,$$ 
where $\bf F$ is a non-Archimedean field for which ${\bf K}_s\subset \bf F$,
$f_h(g):=f(h^{-1}g)$ for each $g\in G$.
For $\mu \in M_l(G,H)$ and $\nu \in M(H)$  let 
$$(\nu * \mu )(A):=\int_H\mu _h(A)\nu (dh) \mbox{ and }
(q\tilde *f)(g):=\int_Hf(hg)q(h)\nu (dh)$$
be convolutions of measures and functions, where $M(H)$ is the space of 
tight measures on $H$ with a finite norm, $\nu \in M(H)$ and
$q\in L(H,\nu ,{\bf F})$.
\par {\bf 2. Lemma.} {\it The convolutions are continuous 
$\bf F$-linear mappings
$$*: M(H)\times M_l(G,H)\to M_l(G,H)\mbox{ and}$$
$$\tilde *: L(H,\nu ,{\bf F})
\times L_H(G,\mu ,{\bf F})\to L_H(G,\mu ,{\bf F}).$$}
\par {\bf Proof.} This follows from Theorem $7.16$ \cite{roo} and estimates
$$\| \nu *\mu \|\le \| \nu \| \times \| \mu \| ,\mbox{ }
\| q\tilde *f\|_{L_H(G,\mu , {\bf F})} \le \| q\| _{L(H,\nu ,{\bf F})}
\times \| f\|_{L_H(G,\mu ,{\bf F})},$$
since $\rho _{\mu }(h,g)\in L(H\times G, \nu \times \mu , {\bf F})$.
\par {\bf 3. Lemma.} {\it For $\mu \in M_l(G,H)$ 
the translation map $(q,f)\to f_q(g)$ is continuous from $H\times 
L_H(G,\mu ,{\bf F})$ into $L_H(G,\mu ,{\bf F})$.}
\par {\bf Proof.} In view of Lemma $7.10$ and Theorem $7.12$ \cite{roo}
for each $\epsilon >0$ the set $\{ x:$ $|f(x)|N_{\mu }(x)\ge \epsilon \} $
is $Af (G,\mu )$-compact and $f$ is $Af (G,\mu )$-continuous.
The embedding of $H$ into $G$ is compact (see \S 1), hence for each
$q\in H$ there exists $V$ clopen in $H$ and such that $q^{-1}V$ is a
subgroup of $H$ with $cl_Gq^{-1}V$ compact in $G$, where $cl_G(A)$
denotes the closure of a subset $A$ in $G$. The product of
compact subsets in $G$ is compact in $G$, hence $f_q(g)\in L_H
(G,\mu ,{\bf F})$ and $(q,f)\mapsto f_q(g)$ is the continuous mapping,
since the restriction of the $Bco (G)$ and the $Af (G,\mu )$-topologies
onto $X_{\epsilon }$ coincide, $\| f_q \|_{L_H(\mu )}= \| f \|_{L_H(\mu )}$
for each $q\in H$ (see \S 1.(b)).
\par {\bf 4. Proposition.} {\it For a probability measure 
$\mu \in M(G)$ there exists
an approximate unit which is a sequence of nonzero continuous 
functions $\psi _i :G \to \bf F$ such that $\int_G\psi _i(g)\mu (dg)=1$
and for each neighbourhood $U\ni e$ in $G$ there exists $i_0$ such that
$supp( \psi _i) \subset U$ for each $i>i_0$.}
\par {\bf Proof.} A group $G$ has a countable base of neighbourhoods
of $e\in G$. A measure $\mu $ is quasi-invariant, hence
$\| U \|_{\mu }>0$ for each neighbourhood $U\ni e$, $\mu $ is the
tight measure, hence there exists a system of neighbourhoods
$ \{ U_i: U_i\ni e \forall i \} $, $\bigcap_iU_i= \{ e \} $,
$U_i\supset U_{i+1}$ for each $i$, $supp (\psi _i)\subset U_i$.
Choose $\psi _i$ such that $\int_G\psi _i(g)\mu (dg)=1$ for each $i$.
\par {\bf 5. Proposition.} {\it If $(\psi _i: i\in {\bf N})$
is an approximate unit in $H$ relative to a probability measure
$\nu \in M(H)$, then $\lim_{i\to \infty }\psi _i*f=f$ in the
$L_H(G,\mu ,{\bf F})$ norm, where $\mu \in M_l(G,H)$, 
$f\in L_H(G,\mu ,{\bf F})$.} 
\par {\bf Proof.} In view of Theorem $7.12$ \cite{roo} for each
$\epsilon >0$ and each $i$ there exists a finite number of $h_j\in H$,
$j=1,...,n,$ $n\in \bf N$, such that $\bigcup_{j=1}^nh_jU_i
\supset X_{\epsilon ,f}\times V,$ where $\{ x:$ $ |f(x)| N_{\mu }(x)
\ge \epsilon \} =:X_{\epsilon ,f}$, $V$ is a clopen neighbourhood of $\xi 
\in H$ which can be chosen such that $\xi ^{-1}V$ is a subgroup of $H$,
$cl_G(\xi ^{-1}V)$ is compact in $G$,
\par $(\psi _i*f)(g)=\int_{h\in U_i}\psi _i(h)f(hg)\nu (dh)$ and \\
$\sup_{g\in G} |f_{\xi }(g)-(\psi _i*f_{\xi })(g)|N_{\mu }(g) \le $
$\sup_{g\in G} [\sup_{h\in X_{\epsilon ,f_{\xi }}} |\psi _i(h)|
|f_{\xi }(hg)-f_{\xi }(g)| N_{\mu }(g) +\epsilon \| \psi _i \|_{\nu }
\| G \|_{\mu }]$, \\
hence $\lim_{i\to \infty } (\psi _i*f)=f$ in $L_H(G,\mu ,{\bf F})$-norm.
\par {\bf 6. Lemma.} {\it Suppose $g\in L_H(G,\mu ,{\bf F})$
and $(g^x|_H)\in L(H,\nu ,{\bf F})$ for each $x\in G$, 
$f\in L(H,\nu ,{\bf F})$, where $g^x(y):=g(yx)$ for each $x$ and $y\in G$. 
Let $\mu $ and $\nu $ be probability measures, $\mu \in M_l(G,H)$,
$\nu \in M(H)$. Then $f{\tilde *}g\in L_H(G,\mu ,{\bf F})$ and there exists
a function $h: G \to \bf F$ such that $h|_H$ is continuous, $h=f{\tilde *}g$
$\mu $-a.e. on $G$ and $h$ vanishes at $\infty $ on $G$.}
\par {\bf Proof.} In view of Fubini theorem we have 
$$\| f{\tilde *}g\|_{L_H(G,\mu ,{\bf F})}=\sup_{h\in H, z\in G}
|\int_Hf(y) g(yz)\nu (dy)| N_{\mu }(z)\le $$
$$\| g(z) \|_{L_H(G,\mu ,{\bf F})} \| f \|_{L(H,\nu ,{\bf F})}$$
The equation $\alpha _f(\phi ):=\int_Hf(y)
\phi (y)\nu (dy)$ defines a continuous linear functional
on $L^{\infty }(H,\nu ,{\bf F}):=\{ \phi : H\to {\bf F}:$ $ \phi $
$\mbox{is}$ $(Af(H,\nu ),Bco({\bf F}))\mbox{-measurable},$
$\| \phi \|_{\infty }:=ess_{\nu }-\sup_{x\in H} |\phi (x)|<\infty \} $.
In view of Lemma 3 the function
$\alpha _f(g^{(qx)^{-1}})=:v(qx)=:w(q,x)$ of two variables
$q$ and $x$ is continuous on $H\times H$ for $q, x \in H$, since the mapping
$(q,x)\mapsto (qx)^{-1}$ is continuous from $H\times H$ into $H$.
By Theorem $7.16$ \cite{roo}
$$\int_Gv(y)\psi (y)\mu (dy)=
\int_G\int_Hf(y)g(yx)\psi (x)\nu (dy) \mu (dx)=$$
$$\int_Hf(y)[\int_Gg(yx) \psi (x)\mu (dx)]\nu (dy)$$ 
for each $\psi \in L^{\infty }(G,\mu ,{\bf F})$.
From this it follows, that $\mu ( \{ y: 
h(y)\ne (f{\tilde *}g)(y), y\in G \} )=0$, since $h$ and $(f{\tilde *}g)$
are $\mu $-measurable functions due to Fubini theorem and 
the continuity of the composition and the inversion in a topological group.
In view of Theorem $7.12$ \cite{roo}
for each $\epsilon >0$ there are compact subsets $C\subset H$ and
$D\subset G$ and functions $f'\in L(H,\nu ,{\bf F})$
and $g'\in L_H(G,\mu ,{\bf F})$ with closed supports
$supp(f')\subset C$, 
$supp(g')\subset D$ such that $cl_GCD$ is compact in $G$,
$$\|f'-f\|_{L(H,\nu ,{\bf F})}
<\epsilon \mbox{ and }\|g'-g\|_{L_H(G,\mu ,{\bf F})} <\epsilon ,$$ 
since by the supposition of \S 1 the group $H$
has the base ${\sf B}_H$ of its topology $\tau _H$, 
such that the closures $cl_GV$ are compact in $G$ for 
each $V\in {\sf B}_H$. From the inequality
$$|h'(x)-h(x)|\le (\| f\|_{L(H,\nu ,{\bf F})}
+\epsilon )\epsilon +\epsilon \| g\|_{L_H(G,\mu ,{\bf F})}$$ it
follows that for each $\delta >0$ there exists a compact 
subset  $K\subset G$ 
with $|h(x)|<\delta $ for each $x\in G\setminus K$, where
$h'(x^{-1}):=\alpha _{f'}({g'}^x).$
\par {\bf 7. Proposition.} {\it Let $A, B\in Af(G,\mu )$, 
$\mu $ and $\nu $ be probability measures, $\mu \in M_l(G,H)$,
$\nu \in M(H)$. Then 
the function $\zeta (x):=\mu (A\cap xB)$ is continuous on $H$
and $\nu (yB^{-1}\cap H)\in L(H,\nu ,{\bf F})$. Moreover, 
if $\| A \|_{\mu } \| B \|_{\mu } >0,$ 
$\mbox{ }\mu ( \{ y \in G: yB^{-1}\cap H \in Af(H,\nu )\mbox{ and }
\| yB^{-1}\cap H \| _{\nu }> 0 \} )>0$, then $\zeta (x)\ne 0$ on $H$.}
\par {\bf Proof.} Let $g_x(y):=Ch_A(y)Ch_B(x^{-1}y)$, then $g_x(y)
\in L_H(G,\mu ,{\bf F})$, where $Ch_A(y)$ is the characteristic function
of $A$. In view of Propositions 4 and 5
there exists $\lim_{i\to \infty }\psi _i*g_x=g_x$ in 
$L_H(G,\mu ,{\bf F})$. In view of Lemma $6$ $\zeta (x)|_H$ is continuous.
There is the following inequality:
$$1\ge |\int_H\mu (A\cap xB)\nu (dx)|=|\int_H\int_GCh_A(y)Ch_B(x^{-1}y)
\mu (dy)\nu (dx)|.$$
In view of Theorem $7.16$ \cite{roo} there exists
$$\int_HCh_B(x^{-1}y)\nu (dy)=\nu ((yB^{-1})\cap H)\in L(G,\mu ,{\bf F}),
\mbox{ hence}$$
$$\int_H\mu (A\cap xB)\nu (dx)=\int_G\nu (yB^{-1}\cap H)
Ch_A(y)\mu (dy).$$
\par {\bf 8. Corollary.} {\it Let $A, B\in Af(G,\mu )$, $\nu \in M(H)$ and 
$\mu \in M_l(G,H)$ be probability measures. Then
denoting $Int_HV$ the interior of a subset $V$ of $H$ 
with respect to $\tau _H$, one has
\par $(i)$ $Int_H(AB)\cap H\ne \emptyset $, when
$$\| \{ y\in G: \| yB\cap H \|_{\nu } >0 \} \|_{\mu } >0;$$
\par $(ii)$ $Int_H(AA^{-1})\ni e$, when 
$$\| \{ y\in G:  \| yA^{-1}\cap H \| _{\nu } >0 \} \| _{\mu } >0.$$}
\par {\bf Proof.} $AB\cap H\supset \{ x\in H: \| (A\cap xB^{-1})
\|_{\mu }>0 \} $.
\par {\bf 9. Corollary.} {\it Let $G=H$. If $\mu \in M_l(G,H)$ 
is a probability measure, then $G$ is a locally compact topological group.}
\par {\bf Proof.} Let us take $\nu =\mu $ and $A=C\cup C^{-1}$, where $C$ 
is a compact subset of $G$ with $\| (C) \| _{\mu }>0$, whence
$\|(yA) \| _{\mu }>0$
for each $y\in G$ and inevitably $Int_G(AA^{-1})\ni e$.
\par {\bf 10. Corollary.} {\it Let $\mu \in M_l(G,H)$,
$\rho _{\mu }(h,z)\in L(H,\nu ,{\bf F})\times L(G,\mu ,{\bf F})$, then
$\rho _{\mu }(h,z)$ is continuous $\nu \times \mu $-a.e. on
$H\times G$.}
\par {\bf Proof.} In view of the cocycle condition \\
$\rho _{\mu }(\phi \psi ,g):=\mu _{\phi \psi }(dg)/\mu (dg)$
$=(\mu _{\phi \psi }(dg)/\mu _{\phi }(dg))(\mu _{\phi }(dg)/\mu (dg))$ \\
$=\rho _{\mu }(\psi , \phi ^{-1}g) \rho _{\mu }(\phi ,g)$ \\
on $\rho _{\mu }$ and Corollary $8$ above, Theorem $7.12$ \cite{roo}
for each $\epsilon >0$ the quasi-invariance factor $\rho _{\mu }(h,g)$
is continuous on $H\times G_{\epsilon }$, but $G_{\epsilon }$
is neighbourhood of $e$ in $G$, where $G_{\epsilon }:=\{ g\in G:
N_{\mu }(g)\ge \epsilon \} $. Since $H$ is dense in $G$ and
$G$ is separable, then
$\bigcup_{j=1}^{\infty }h_jG_{\epsilon }=G$, where $ \{ h_j: $
$j\in {\bf N} \} $ is a countable subset in $G$. Therefore,
$\rho _{\mu }(h,z)$ is continuous $\nu \times \mu $-a.e. on $H\times G$.
\par {\bf 11. Corollary.} {\it Let $G$ be a locally compact group,
$\rho _{\mu }(h,z)\in L(G\times G,\mu \times \mu ,{\bf F})$,
then $\rho _{\mu }(h,z)$ is $\mu \times \mu $-a.e. continuous
on $G\times G$.}
\par {\bf Proof.} It follows from Corollaries $9$ and $10$.
\par {\bf 12. Remark.} The latter two corollaries show, that
the condition of continuity of $\rho _{\mu }(h,z)$ imposed in
\cite{luumn582,lulapm} is not restrictive.
\par {\bf 13. Lemma.} {\it Let $\mu \in M_l(G,H)$ be a probability measure
and $G$ be non-locally compact. Then $\| H \| _{\mu }=0$.}
\par {\bf Proof.} This follows from Theorem 3.13 \cite{lulapm} above and
the proof of Lemma 2, since the embedding $T_eH\hookrightarrow T_eG$
is a compact operator in the non-Archimedean case and a tight
measure $\mu $ on $G$ induces a tight measure on a neighbourhood
$V$ of $0$ in $T_eG$ such that $V$ is topologically homeomorphic
to a clopen subgroup $U$ in $G$ (see also papers about
construction of quasi-invariant measures on the considered here groups
\cite{lu1,lutmf99,luumn983,lubp2,luseamb2}).
\par  {\bf 14. } Let $(G,\tau _G)$ and $(H,\tau _H)$ be a pair of
topological non-locally compact groups $G,$ $H$ (Banach-Lie, Frechet-Lie
or groups of diffeomorphisms or loop groups) with uniformities
$\tau _G, \tau _H$ such that $H$ is dense in $(G,\tau _G)$ and there is
a probability measure $\mu \in M_l(G,H)$ with continuous
$\rho _{\mu }(z,g)$ on $H\times G$.
Also let $X$ be a Banach space over $\bf F$ and $IS (X)$ be
the group of isometric $\bf F$-linear automorphisms of $X$
in the topology inherited from the Banach space $L(X)$
of all bounded $\bf F$ linear operators from $X$ into $X$.
\par {\bf Theorem.} {\it $(1).$ If $T: G\to IS(X)$ is a weakly continuous
representation, then there exists $T': G\to IS(X)$ equal $\mu $-a.e. to
$T$ and $T'|_{(H,\tau _H)}$ is strongly continuous.
\par  $(2).$ If $T: G\to IS(X)$ is a weakly measurable representation 
and $X$ is of separable type $c_0 ({\bf F})$ over $\bf F$,
then there exists $T': G\to IS(X)$ equal to $T$ $\mu $-a.e. and
$T'|_{(H,\tau _H)}$ is strongly continuous.}
\par  {\bf Proof.} Take ${\cal K}(G):=(I)\cup L(G,\mu ,{\bf F})$,
where $I$ is the unit operator on $L(G,\mu ,{\bf F})$. Then we can define 
$$A_{(\lambda e+a)_h}:=
\lambda I+ \int_Ga_h(g)[\rho _{\mu }(h,g)]T_g\mu (dg),$$ 
where $a_h(g):=a(h^{-1}g)$. Then 
$$|(A_{(\lambda e+a)_h}-A_{\lambda e+a}\xi ,
\eta )|\le \sup_{g\in G} [|a_h(g)\rho _{\mu }(h,g)-a(g)|
\mbox{ }|T_g\xi |N_{\mu }(g)],$$ 
hence $A_{a_h}$ is strongly continuous with respect to 
$h\in H$, that is, 
$$\lim_{h\to e}|A_{a_h}\xi -A_a\xi |=0.$$
Denote $A_{a_h}=T'_hA_a$, so $T'_h\xi =A_{a_h}
\xi _0$, where $\xi =A_a\xi _0$, $a\in L(G,\mu ,{\bf F})$. Whence 
$${T'}_h\xi=A_{a_h}\xi _0=\int _Ga_h(g)T_g\xi _0\rho _{\mu }(h,g)\mu (dg)$$ 
$$= T_h\int _G a(z) T_z\xi _0\mu (dz)=T_h\xi ,$$
hence $|{T'}_h\xi |=|\xi |$ for each $h\in H$. Therefore,
$T'_h$ is uniquely extended to an isometric operator on the Banach space
$X' \subset X$. In view of Lemma 10, $\| H \|_{\mu }=0$.
Hence $T'$ may be considered
equal to $T$ $\mu $-a.e. Then a space $sp_{\bf F}[A_{a_h}: h\in H]$
is evidently dense in $X$, since 
$$A_{a_h}\xi =\int_Ga_h(g)T_g\rho _{\mu }(h,g)\mu (dg)\xi =$$
$$T_h\int_G a(g)T_g\mu (dg)\xi .$$
For proving the second statement let \\
${\sf R}:=[\xi : A_a\xi =0\mbox{ for each }a\in L(G,\mu ,{\bf F})]$.
If $\eta \in X^*$ and
$$\eta (A_a\xi )=\eta (\int_Ga(g)T_g\xi \mu (dg))=
\eta (\int_Ga(g){T'}_g\xi \mu (dg))$$ 
for each $a(g)\in L(G,\mu ,{\bf F})$, then $\eta (T_g\xi )=\eta (T'_g\xi )$
for $\mu $-almost all $g\in G$, where $X^*$ denotes the topological
dual space of all $\bf K$-linear functionals $f: X\to \bf K$.
Suppose that $\{ \xi _n^*:$ $n\in {\bf N} \} $ is an
orthonormal system in $X^*$ separating points of $X$.
It exists, since by the supposition of this theorem
$X=c_0({\bf F})$. If $\xi \in X$, then
$$\xi _m^* (\int_Ga(g)T_g\xi  \mu (dg))=0$$ 
for each $g\in G\setminus S_m$, where $\| S_m \|_{\mu }=0$. 
Therefore, $\xi _m^*(T_g\xi )=0$ for each $m\in \bf N$, if 
$g\in G\setminus S$, where $ S:=\bigcup_{m=1}^{\infty }S_m$.
Hence $T_g\xi =0$ for each $g\in G\setminus S$, consequently,
$\xi =0$. Consider the embedding $X\hookrightarrow X^*$
with the help of the standard orthnormal basis $\{ e_j: j \} $ in
$X$ over $\bf F$. Therefore, let $\{ \xi _m: m\in {\bf N} \}
\subset X\hookrightarrow X^*.$
Then $\xi _m^*(T_g\xi _n)=\xi _m^*(T'_g\xi _n)$
for each $g\in G\setminus \gamma _{n,m}$, where
$ \| \gamma _{n,m} \|_{\mu }=0$, $\xi _m^*$ is the image of $\xi _m$
under this standard embedding $X\hookrightarrow X^*$.
Hence $\xi _m^*(T_g\xi _n)=\xi _m^*(T'_g\xi _n)$ for each $n, m\in \bf N$
and each $g\in G\setminus \gamma $, where $\gamma :=
\bigcup_{n,m}\gamma _{n,m}$ and $\| \gamma \|_{\mu }=0$,
since in $L(G,\mu ,{\bf F})$ the family of all step functions is dense.
Therefeore, ${\sf R}=0$.
\par {\bf 15. Definition and note.} Let $\{ G_i: i\in {\bf N_o} \}$
be a sequence of topological groups such that $G=G_0$,
$G_{i+1}\subset G_i$ and
$G_{i+1}$ is dense in $G_i$ for each $i\in \bf N_o$
and their topologies are denoted $\tau _i$, $\tau _i|_{G_{i+1}}
\subset \tau _{i+1}$ for each $i$, where $N_o:=\{ 0,1,2,... \} $.
Suppose that these groups are supplied with $\bf F$-valued probability
quasi-invariant measures
$\mu ^i$ on $G_i$ relative to $G_{i+1}$. For example, such sequences exist
for groups of diffeomorphisms or loop groups considered in previous papers
\cite{lu1,lutmf99,luumn983,lubp2,luseamb2}).
Let $L_{G_{i+1}}(G_i,\mu ^i,{\bf F})$ denotes a Banach subspace
of $L(G_i,\mu ^i,{\bf F})$ as in \S 1(b). Let ${\tilde L}(G_{i+1},
\mu ^{i+1},L(G_i,\mu ^i,{\bf F}))=: H_i$ denotes the completion of the
subspace of $L(G_i,\mu ^i,{\bf F})$ of all elements $f$ such that
$$\| f\|_i:=\max [{\| f^2\|^{1/2}}_{L(G_i,\mu ^i,{\bf F})},
{\| f\| '}_i]<\infty ,\mbox{ where}$$
$${\| f\| '}_i :=[\sup_{x\in G_i, y\in G_{i+1}}|f(y^{-1}x)|^2
N_{\mu ^i}(x)\max (1, N_{\mu ^{i+1}}(y)) ]^{1/2}.$$ 
Evidently $H_i$ are Banach spaces over $\bf F$. Let 
$$f^{i+1}*f^i(x):=\int_{G_{i+1}}f^{i+1}(y)f^i(y^{-1}x)\mu ^{i+1}(dy)$$
denotes the convolution of $f^i\in H_i$.
\par {\bf 16. Lemma.} {\it The convolution $*: H_{i+1}\times H_i\to H_i$
is the continuous $\bf F$-bilinear mapping.}
\par {\bf Proof.} From the definitions we have: \\
$\| (f^{i+1}*f^i)^2 \|^{1/2}_{L(G_i,\mu ^i,{\bf F})}=
\sup_{x\in G_i} |(f^{i+1}*f^i)(x)| N^{1/2}_{\mu ^i}(x)$ \\
$\le \sup_{x\in G_i, y\in G_{i+1}}[|f^{i+1}(y)|N^{1/2}_{\mu ^{i+1}}(y)]
[|f^i(y^{-1}x)|N^{1/2}_{\mu ^i}(x)N^{1/2}_{\mu ^{i+1}}(y)]$ \\
$=\| (f^{i+1})^2 \|^{1/2}_{L(G_{i+1},\mu ^{i+1},{\bf F})}
{ \| f^i \| '}_i$ and \\
${ \| f^{i+1}*f^i \| '}_i =[\sup_{x\in G_i, y\in G_{i+1}}
|(f^{i+1}*f^i)(y^{-1}x)|^2N_{\mu ^i}(x) N_{\mu ^{i+1}}(y)]^{1/2}$ \\
$\le [\sup_{x\in G_i, y\in G_{i+1}, z\in G_{i+1}} |f^{i+1}(z)|^2
N_{\mu ^{i+1}} (z) |f^i(y^{-1}z^{-1}x)|^2 N_{\mu ^i}(x)
N_{\mu ^{i+1}}(y)N_{\mu ^{i+1}}(z)]^{1/2}$ \\
$\le \| (f^{i+1})^2 \| ^{1/2}_{L(G_{i+1},\mu ^{i+1},{\bf F})}
{ \| f^i \| '}_i,$ \\
since from $A\ni y$ and $B\ni z$ for $A, B\in Af (G_{i+1},\mu ^{i+1})$
it follows that $AB\ni yz$ and $AB\in Af (G_{i+1},\mu ^{i+1})$,
which follows from \\
$\mu ^{i+1}(AB)=\int_{A\ni a} \mu ^{i+1}(aB)\mu ^{i+1}(da)=
\int_{B\ni b}\int_{A\ni a} \mu ^{i+1}(adb)\mu ^{i+1}(da),$ so that \\
$\| A \|_{\mu ^{i+1}} \le \| G_{i+1} \| _{\mu ^{i+1}}=1$, hence
$N_{\mu ^{i+1}}(z)\le 1$ for each $z\in G_{i+1}$
for the probability measure $\mu ^{i+1}$. Therefore,
$\| f^{i+1}*f^i \| _i\le \| f^{i+1} \|_{i+1} \| f^i \| _i$.
\par {\bf 17. Definition.} Let $c_0(\{ H_i: i\in {\bf N_o} \} )=:H$
be the Banach space consisting of elements
$f=( f^i: f^i\in H_i, i\in {\bf N_o})$, 
for which $ \lim_{i\to \infty } \| f^i \|_i=0$, where
$$\| f\|:=\sup_{i=0}^{\infty }\| f^i\|_i<\infty .$$
For elements $f$ and $g\in H$ their convolution is defined 
by the formula: $f\star g:=h$ 
with $h^i:=f^{i+1}*g^i$ for each $i\in \bf N_o$.
Let $*: H\to H$ be an involution such that $f^*:=({f^j}^{\wedge }: 
j\in {\bf N_o} )$, where ${f^j}^{\wedge }(y_j):=f^j(y_j^{-1})$ for each
$y_j\in G_j$, $f:=(f^j: j\in {\bf N_o} )$.
\par {\bf 18. Lemma.} {\it $H$ is a non-associative non-commutative
Banach algebra with involution $*$, that is $*$ is $\bf F$-bilinear
and $f^{**}=f$ for each $f\in H$.}
\par {\bf Proof.} In view of Lemma 16 the convolution
$h=f\star g$ in the Banach space $H$
has the norm $\| h\| \le \| f\|$ $\| g\| $, hence
is a continuous mapping from $H\times H$ into $H$.
From its definition it follows that the convolution is
$\bf F$-bilinear.
It is non-associative as follows from the computation of
i-th terms of $(f\star g)\star q$ and $f\star (g 
\star q)$, which are
$(f^{i+2}*g^{i+1})*q^i$ and $f^{i+1}*(g^{i+1}*q^i)$ respectively,
where $f$, $g$ and $q\in H$.
It is non-commutative, since there are $f$ and $g\in H$ for which
$f^{i+1}*g^i$ are not equal to $g^{i+1}*f^i$.
From ${f^j}^{\wedge \wedge }(y_j)=f^j(y_j)$ it follows that
$f^{**}=(f^*)^*=f$. 
\par {\bf 19. Note.} In general $(f\star g^*)^*\ne g\star f^*$
for $f$ and $g\in H$, since there exist $f^j$ and $g^j$
such that $g^{j+1}*(f^j)^*\ne (f^{j+1}*(g^j)^*)^*$.
If $f\in H$ is such that $f^j|_{G_{j+1}}=f^{j+1}$, then
$$((f^{j+1})^**f^j)(e)=\int_{G_{j+1}}(f^{j+1}(y))^2
\mu ^{j+1}(dy), \mbox{ hence } $$
$$ |((f^{j+1})^**f^j)(e)| \le \| (f^{j+1})^2
\|^{1/2}_{L(G_{j+1},\mu ^{j+1},{\bf F})} \le \| f^{j+1} \|_{j+1},$$
where $j\in \bf N_o$.
\par {\bf 20. Definition.} Consider the standard Banach space
$c_0({\bf F})$ over the field $\bf F$ as a Banach algebra with
the convolution $\alpha \star \beta =\gamma $
such that $\gamma ^i:=\alpha ^{i+1}\beta ^i$, where
$\alpha :=(\alpha ^i: \alpha ^i\in {\bf F}, i\in {\bf N_o})$,
$\alpha $, $\beta $ and $\gamma \in c_0({\bf F})$. 
\par {\bf 21. Note.} The algebra $c_0({\bf F})$ has two-sided ideals
$J_i:=\{ \alpha \in c_0({\bf F}): \alpha ^j=0$ for each $j>i \} $,
where $i\in \bf N_o$. That is, $J\star c_0({\bf F})\subset J$ 
and $c_0({\bf F}) \star J=J$ and $J$ 
is the $\bf F$-linear subspace of $c_0({\bf F})$, but
$J\star c_0({\bf F})\ne J$. 
There are also right ideals, which are not left ideals:
$K_i:=\{ \alpha \in c_0({\bf F}): \alpha ^j=0$ 
for each $j=0,...,i \} $, where
$j\in \bf N_o$. That is, $c_0({\bf F})\star 
K_i=K_i$, but $K_i\star c_0({\bf F})=
K_{i-1}$ for each $i\in \bf N_o$, where $K_{-1}:=c_0({\bf F})$.
The algebra $c_0({\bf F})$ is the particular case of
$H$, when $G_j=\{ e \}$ for each $j\in \bf N_o$.
We consider further $H$ for non-trivial topological groups outlined above
with $G_{\infty }:=\bigcap_{j=0}^{\infty } G_j$ dense in each $G_j$.
\par {\bf 22. Theorem.} {\it If $\cal F$ is a maximal proper left or right
ideal in $H$,
then $H/\cal F$ is isomorphic as the nonassociative noncommutative algebra 
over $\bf F$ with $c_0({\bf F})$.}
\par {\bf Proof.} The ideal $\cal F$ is also the $\bf F$-linear subspace
of $H$. In view of Theorem $7.12$ \cite{roo} a function $f: G_j\to \bf F$
is $\mu ^j$-integrable if and only if it satifies two properties:
$f$ is $Af (G,\mu ^j)$-continuous and for each $\epsilon >0$
the set $ \{ x: |f(x)| N_{\mu ^j}(x) \ge \epsilon \} $
is $Af (G,\mu ^j)$-compact and hence contained in 
$\{ x: N_{\mu ^j} (x)\ge \delta \} $ for some $\delta >0$.
Suppose, that there exists $j\in \bf N_o$ such that $f^j=0$ for each
$f\in \cal F$, then $f^i=0$ for each $i\in \bf N_o$,
since the space of bounded
$\bf F$-valued continuous functions $C^0_b(G_{\infty }, {\bf F})$
on $G_{\infty }$ is dense 
$H_j:=\{ f^j: f\in H \} $ and $C^0_b(G_{\infty },{\bf F})\cap F_j=\{ 0\} $
and $C^0_b(G_j,{\bf F})|_{G_{j+1}}\supset C^0_b(G_{j+1}, {\bf F}).$
Therefore, ${\cal F}_j\ne \{ 0 \} $ for each $j\in \bf N_o$,
where ${\cal F}_j:= \{ f^j: f\in {\cal F} \} $, consequently,
${\bf F}\hookrightarrow F_j$ for each $j\in \bf N_o$.
Since ${\bf F}$ is embeddable into each $F_j$, then
there exists the embedding of $c_0({\bf F})$ into $\cal F$,
where $H_j:=\{ f^j: f\in H \} $, $\pi _j: H\to H_j$ are 
the natural projections.
\par  The subalgebra $\cal F$ is closed in $H$, since
$H$ is the topological algebra and $\cal F$ is the maximal proper subalgebra.
The space $H_{\infty }:=\bigcap_{j\in \bf N_o}H_j$
is dense in each $H_j$. 
\par If ${\cal F}_i=H_i$ for some $i\in \bf N_o$,
then ${\cal F}_j=H_j$ for each $j\in \bf N_o$, since
$C^0_b(G_{\infty },{\bf F})$ is dense in each $H_j$
and $C^0_b(G_j,{\bf F})|_{G_{j+1}}\supset C^0_b(G_{j+1}, {\bf F}).$
The ideal $\cal F$ is proper, consequently, ${\cal F}_j\ne H_j$
as the $\bf F$-linear subspace for each $j\in \bf N_o$,
where ${\cal F}_j=\pi _j({\cal F})$. 
\par There exist $\bf F$-linear continuous operators from
$c_0({\bf F})$ into $c_0({\bf F})$ such that
$x\mapsto (0,...,0,x^0,x^1,x^2,...)$
with $0$ as $n$ coordinates at the beginning,
$x\mapsto (x^n,x^{n+1},x^{n+2},...)$ for $n\in \bf N$;
$x\mapsto (x^{kl+\sigma _k(i)}: k\in {\bf N_o}, i\in (0,1,...,l-1) )$,
where ${\bf N}\ni l\ge 2$, $\sigma _k\in S_l$ are elements of the symmetric
group $S_l$ of the set $(0,1,...,l-1)$.
Then $f\star (g\star h)+c_0({\bf F})$ and $(f\star g)\star h+c_0({\bf F})$ 
are considered as the same class, also $f\star g+c_0({\bf F})=
g\star f+c_0({\bf F})$ in $H/c_0({\bf F})$,
since $(f+c_0({\bf F}))\star (g+c_0({\bf F}))=f\star g+c_0({\bf F})$
for each $f, g$ and $h\in H$. Then $f\star (g\star h)+c_0({\bf F})$
and $(f\star g)\star h+c_0({\bf F})$ are considered as the same class
for each $f, g, h\in \cal F$, also $f\star g+c_0({\bf F})=g\star f+
c_0({\bf F})$ in ${\cal F} / c_0({\bf F})$,
since $(f+c_0({\bf F}))\star (g+c_0({\bf F}))=f\star g
+l_2({\bf F})\subset \cal F$ for each $f$ and $g\in \cal F$. 
Therefore, the quotient algebras $H/c_0({\bf F})$
and ${\cal F}/ c_0 ({\bf F})$ are the associative commutative
Banach algebras.
\par From $\mu ^i\in M_l(G_i,G_{i+1})$ it follows that for each open
subset $W\ni e$, $W\subset G_i$ there exists a clopen subgroup $U\subset
W$ such that $\mu ^i(U)\ne 0$, since otherwise $\mu ^i(zV)=0$ for
each $z\in G_{i+1}$ and each open $V\subset W$, hence $\mu ^i(G_i)=0$
contradicting supposition, that each $\mu ^i$ is the probability measure.
\par Let us adjoin a unit to $H/c_0 ({\bf F})$ and to
${\cal F}/c_0 ({\bf F})$.
There is satisfied the equality $Ch_{G_{i+1}}*Ch_{G_i}=Ch_{G_i}$.
Let $U_{i,j}$ be a clopen subgroup in $G_i$, that is possible,
since each $G_i$ is ultrametrizable. Choose $U_{i,j}$ such that
$U_{i,j+1}\subset U_{i,j}\cap G_{j+1}$ and $U_{i,j}\supset U_{i+1,j}$
for each $i$ and $j$, $\bigcap_iU_{i,j}=e\in G_j$ for each $j$.
Since $\mu ^j(G_j)=1$ and $\| G_j \| _{\mu ^j}=1$, then
by induction $(i,j)\in \{ (1,1), (1,2), ..., (1,n),...;
(2,1), (2,2),..., (2,n),...; ... (m,1), (m,2),..., (m,n),... \} $,
where $m, n\in \bf N$, there exists a family $\alpha _{i,j}\in \bf F$
and $ \{ U_{i,j}: i, j \} $ such that
$\alpha _{i,{j+1}}Ch_{U_{i,{j+1}}}*\alpha _{i,j}Ch_{U_{i,j}}=
\alpha _{i,j}Ch_{U_{i,j}}$ and $0< |\alpha _{i,j}| |\mu ^i(U_{i,j})|\le 1$
for each $i, j$.
Put $e_i:= \{ \alpha _{i,j}Ch_{U_{i,j}}: j\in {\bf N_o} \} $,
then $e_i*e_i=e_i$ for each $i$. From the properties of $U_{i,j}$
it follows, that $span_{\bf F} \{ e_i(z^{-1}g): i\in {\bf N},
z\in G_{\infty } \} $ is dense in $H$, where $g:= (g_j: g_j\in G_j
\forall j \in {\bf N_o} )$, $z^{-1}g=(z^{-1}g_j: j)$.
\par Consider the algebras $H/c_0({\bf F})=:A$ and ${\cal F}/c_0({\bf F})=:B$.
The algebras $A$ and $B$ are commutative and associative.
From the preceding proof it follows that $span_{\bf F} \{ e_i(z^{-1}g)+
c_0({\bf F}):$ $i\in {\bf N},$ $z\in G_{\infty } \} $ is dense in $A$,
each $e_i(z^{-1}g)+c_0({\bf F})$ is the idempotent element in $A$.
Therefore, by Theorem $6.12$ \cite{roo} $A$ is the $C$-algebra.
By the definition this means, that there exists a locally compact
zero-dimensional Hausdorff space $X$ such that $A$ is isomorphic with
$C_{\infty }(X,{\bf F})$, where $C_{\infty }(X,{\bf F})$
is the subspace of all $f\in C_b(X,{\bf F})$ for which
for each $\epsilon >0$ there exists a compact subset $X_{\epsilon ,f}$
of $X$ with $|f(x)|<\epsilon $ for each $x\in X\setminus X_{\epsilon ,f}$.
In accordance with Theorem $6.3$ \cite{roo} each maximal
ideal $\cal B$ of $C_{\infty }(X,{\bf F})$ has the form
${\cal B}=\{ f\in C_{\infty }(X,{\bf F}):$ $f(z_0)=0 \} ,$
where $z_0$ is a marked point in $X$.
On the other hand, as it was proved above
${\cal F}_j\ne H_j$ for each $j\in \bf N_o$,
hence there exists the following embedding
$c_0({\bf F})\hookrightarrow (H/{\cal F})$ and $(H/{\cal F})/
c_0({\bf F})$ is isomorphic with $(H/c_0({\bf F}))/({\cal F}/c_0({\bf F}))$.
Therefore, $H/\cal F$ is isomorphic with $c_0({\bf F})$.
\par {\bf 23.} {\large Comments.} Another methods of construction
of isometrical representations of topological totally disconnected
groups which may be nonlocally compact with the help of quasi-invariant
$\bf F$-valued measures were given in
\cite{luseamb,luseamb2,lubp2,lutmf99,luumn01,luijmms3,lufpmsp,luddiagr}
and references therein.
\par Using results of this article it is possible
to make further investigations of nonlocally compact totally
disconnected topological groups, their structures and representations,
measurable operators in Banach spaces over non-Archimedean fields,
apply this for the development of non-Archimedean quantum mechanics
and quantum field theory, quantum gravity, superstring theory
and gauge theory, etc.

\end{document}